\documentclass{article}
\usepackage{amssymb}

\newcommand{\Z}{{\mathbb Z}}
\newcommand{\C}{{\mathbb C}}

\newcommand{\R}{{\mathbb R}}

\newcommand{\La}{\Lambda}

\newcommand{\Ss}{\mathcal{S}} \newcommand{\B}{\mathcal{B}}\newcommand{\K}{\mathcal{K}}

\begin{document}

\begin{large}
\title{\bf On Beurling's sampling theorem in $\R^n$}

\author{Alexander Olevskii  and Alexander
Ulanovskii}

\maketitle

\begin{abstract}
We present an elementary proof of the classical Beurling sampling theorem  which gives a  sufficient condition for sampling of multi--dimensional band--limited functions.
\end{abstract}


\section{Introduction}

 Let $\Ss\subset \R^n, n\geq1,$ be a compact. The  Bernstein space  $ B_\Ss$ consists of all
 bounded functions  on $\R^n$ whose spectrum belongs to $\Ss$. The latter means that
$$
\int_{\R^n}f(x)\hat\varphi(x)\,dx=0, \ \ f\in B_\Ss,
$$
for every smooth function $\varphi(x)$ whose support belongs to a ball disjoint from $\Ss.$
Here $\hat \varphi$  denotes the Fourier transform 
$$
\hat\varphi(x)=\int_{\R^n}e^{-i t\cdot x}\varphi(t)\,dt.
$$

A set
$\La\subset \R^n$ is called a sampling set for $B_S$, if there is a positive constant $C$ such that
$$
\Vert f\Vert_\infty\leq C\Vert f|_\La\Vert_\infty, \ \mbox{for every } \ f\in B_S,
$$
where
$$
\Vert f\Vert_\infty:=\sup_{x\in \R^n}|f(x)|, \ \ \Vert f|_\La\Vert_\infty:=\sup_{\lambda\in\La}|f(\lambda)|.
$$
It is a classical problem  to determine when  $\La$ constitutes a sampling set for $B_\Ss.$  Beurling  discovered the importance of the 
{\it lower uniform density $D^-(\La)$} of $\La$ for this problem:
$$D^-(\La):= \lim_{r\to\infty}\frac{\min_{x\in\R^n}\mbox{Card}(\La\cap (x+r\B))}{|r\B|},
$$where $\B$ is  the unit ball in $\R^n$, $x+r\B$ is the ball of radius $r$ centered at $x$  and $|\Ss|$ denotes the measure of a set $\Ss$. In  \cite{b1} he proved the following

\medskip\noindent
{\bf Theorem 1} {\sl Let $\Ss=[a,b]\subset\R$. Then $\La\subset\R$ is a sampling set for $B_{\Ss}$ if and only if
 \begin{equation}D^-(\La)>|\Ss|/{2 \pi}.\end{equation}}

Hence, when $\Ss$ is an interval in $\R$,  the sampling problem can be solved in terms of the density  $D^-(\La).$ Condition $D^-(\La)\geq |\Ss|/(2 \pi)^n$ remains {\it necessary} for sampling in $B_\Ss$, for every compact set $\Ss\subset\R^n.$
This follows from a general result of Landau \cite{l}. On the other hand, simple examples show that in  dimension one condition (1) ceases to be {\it sufficient} already when $\Ss$ is a union of two intervals. 

A new phenomenon occurs in several dimensions: Even for the simplest sets $\Ss$ like a ball or a cube,
    no sufficient conditions for sampling in $B_\Ss$ can be expressed in terms of  $D^-(\La)$.
 The reason for that is that 
the zeros of the multi-dimensional entire functions are not discrete.  One can check that if  $\Ss\subset\R^n$  contains at least two points, then there are functions $f\in B_\Ss$ whose zero set contains sets $\La\subset\R^n$ with arbitrarily large $D^-(\La)$. Clearly, if a function $f\in B_\Ss$ vanishes on $\La$, then $\La$ is not a sampling set for $B_\Ss$
(see also discussion  in \cite{s}, pp. 122--123).

In \cite{b} Beurling obtained the following sufficient condition for sampling in  $B_\B$:

\medskip\noindent
{\bf  Theorem 2} 
    {\sl Assume  $\La\subset\R^n,n\geq1,$ and  $\rho<\frac{\pi}{2}$ satisfy
$$
\La+\rho\B=\R^n.
$$  Then
\begin{equation}
        \Vert f\Vert_\infty\leq \frac{1}{1-\sin \rho} \Vert f|_\La\Vert_\infty,\ \mbox{ for every } f\in B_{\B},
\end{equation}
     and so   $\La$ is a sampling set for  $B_\B$. }

     \medskip 
    In fact,  Beurling in \cite{b} proves a result on balayage of Fourier--Stieltjes transforms which is equivalent to Theorem 2: {\it For every Dirac's measure $\delta_\xi$, there exists a finite measure with masses on $\La$ such that the values of their  Fourier--Stieltjes transforms  agree in the ball $\B$}. We use a completely different elementary approach which allows us to get a more general result, see Theorem 3 below. We shall see that unlike the case of interpolation in several dimensions (see \cite{o}), the "Beurling-type" sampling  is in fact a one-dimensional phenomenon.

   Observe that condition $\La+\rho\B=\R^n$ in Theorem 2 means that  $\La$ is an $\rho$-net, i.e. for every $x\in\R^n$ there exists $\lambda\in\La$ with $|x-\lambda|\leq\rho.$
Hence,   every $\rho$-net with $\rho<\pi/2$ is a sampling set for $B_\B$. This is sharp: Beurling shows that the theorem ceases to be true for $\pi/2-$nets.

Let us in what follows denote by $\K$  a closed convex central-symmetric body with positive measure.   Then
       $$
       \K^o:=\{x\in \R^n: x\cdot t\leq 1 \mbox{ for all } t\in \K\}
       $$
       denotes the polar body of $\K$. In particular, we have $\B^o=\B.$

The following propositions are  formulated in \cite{b} without  proof: 
 
 (i) Estimate (2) in Theorem 2 can be replaced with a better one:
\begin{equation}\label{bb}
     \Vert f\Vert_\infty\leq \frac{1}{\cos \rho}\Vert f|_\La\Vert_\infty.
     \end{equation}

 (ii) Every set $\La$ satisfying $\La+\rho\K^o=\R^n$ with some $\rho<\pi/2$ is a sampling set for $B_\K.$

 We show that estimate (3) holds for every convex central-symmetric body $\K$:
 
\medskip\noindent
      {\bf  Theorem 3} {\sl Assume $\La\subset\R^n$ and  $\rho<\frac{\pi}{2}$ satisfy
\begin{equation}\label{b}
    \La+\rho\K^o=\R^n.
\end{equation}       Then (\ref{bb}) is true,  and so
       $\La$ is a sampling set for $B_\K$. }

\medskip

Clearly, condition (\ref{b}) means that for every $x\in\R^n$ there exists $\lambda\in\La$ such that $\Vert x-\lambda\Vert_{\K^o}\leq\rho$, where $\Vert x\Vert_{\K^o}:=\inf_{a>0}\{x\in a \K^o\}$. Hence, every 
$\rho$-net in the  norm $\Vert\cdot\Vert_{\K^o}$ is a sampling set for $B_\K$ provided $\rho<\pi/2$. 
This  is sharp:

\medskip

\noindent {\bf Proposition 1}. {\sl Suppose a closed convex central-symmetric body  $\Ss$ contains a point $x_0$ with $\Vert x_0\Vert_{\K^o}=\pi/2$. Then there exists $\La\subset\R^n$ with  $\La+\Ss=\R^n$ and a function $f\in B_\K$ such that $f(\lambda)=0,\lambda\in\La$.}

\medskip 
\noindent {\bf Corollary 1}. 
{\sl Suppose a closed convex central-symmetric body  $\Ss$  has the property that every set $\La\subset\R^n$ satisfying  $\La+\Ss=\R^n$ is a sampling set for $B_\K$. Then  $\Ss\subset \rho \K^o$ for some  $\rho<\pi/2.$}

   \section{Proofs}

\noindent{\bf 1. Proof of Proposition 1}.
 By assumption, there exist  $x_0\in\Ss$ and $t_0\in \K$ such that $x_0\cdot t_0=\pi/2.$ The spectrum of the function $\sin(x\cdot t_0)$ consists of two points $\pm t_0\in\K,$ and so $\sin(x\cdot t_0)\in B_\K.$

Let $\La:=\{x\in \R^n: x\cdot t_0\in \pi\Z\}$ be the zero set of  $\sin(x\cdot t_0)$. Denote by
 $I=\{\tau x_0: -1\leq\tau\leq 1\}\subseteq \Ss$ the interval from $-x_0$ to $x_0$. Clearly, for every point $y\in\R^n$ there exist $n\in\Z$ and $-1\leq\tau\leq 1$ such that
$y\cdot t_0 =\pi n-\tau \pi/2.$ Hence, $y-\tau x_0\in \La$, which implies $\La+I=\R^n$. $\Box$

  \medskip\noindent{\bf 2. Proof of Theorem 3}. 
We shall deduce Theorem 3  from the following 

\medskip\noindent
{\bf Lemma 1} {\it Suppose a function $g\in B_{[-\tau,\tau]}$ satisfies $|g(0)|=\Vert g\Vert_\infty$. Then }
\begin{equation}|g(u)|\geq|g(0)|\cos(\tau u), \ \  |u|<\pi/2\tau.\end{equation} 


This lemma is proved in [3, proof of Theorem 4]. For completeness of presentation, we  sketch the proof below.

Let us now prove Theorem 3. Take any function $f\in B_\K$. Assume first that $|f|$ attains maximum on $\R^n$, i.e. $|f(x_0)|=\Vert f\Vert_\infty$ for some $x_0\in\R^n$.  By (4),  there exists $\lambda_0\in\La$ with $\Vert\lambda_0-x_0\Vert_{\K^o}\leq\rho$. 
  Consider the function of one variable $g(u):=f(x_0+u(\lambda_0-x_0)), u\in\R$.
 One may check that $g\in B_{[-\tau,\tau]}$ with $\tau=\Vert\lambda_0-x_0\Vert_{\K^o}$. Also, clearly $|g(0)|=\Vert g\Vert_\infty$ and $g(1)=f(\lambda_0)$. Since $\tau\leq\rho<\pi/2$, we may use inequality (5) with $u=1$:
 $$\Vert f\Vert_\infty=|f(x_0)|=|g(0)|\leq  \frac{|g(1)|}{\cos \tau}\leq
 \frac{|f(\lambda_0)|}{\cos \rho}\leq\frac{1}{\cos \rho}\Vert f|_\La\Vert_\infty. $$

If $|f|$ does not attain maximum on $\R^n$, we consider the function $f_\epsilon(x):=f(x)\varphi(\epsilon x)$, where
$\varphi\in B_{\epsilon\B}$ is any function satisfying $\varphi(0)=1$ and $\varphi(x)\to0$ as $|x|\to\infty.$
It is clear that $f_\epsilon\in B_{\K+\epsilon\B}$ and that $f_\epsilon$ attains maximum on $\R^n$. Set $g_\epsilon(u):=f_\epsilon(x_0+u(\lambda_0-x_0)), u\in\R,$ where $x_0$ and $\lambda_0$ are chosen so that $|g_\epsilon(0)|=\Vert f_\epsilon\Vert_\infty$ and  $\Vert\lambda_0-x_0\Vert_{\K^o}\leq\rho$. We have $g\in B_{[-\tau-\delta,\tau+\delta]} $,
 where $\tau=\Vert\lambda_0-x_0\Vert_{\K^o}\leq\rho<\pi/2$ and $\delta=\delta(\epsilon)\to0$ as $\epsilon\to0.$ So, if $\epsilon$ is so small that $\tau+\epsilon<\pi/2$, we may repeat the argument above to 
obtain $\Vert f_\epsilon\Vert_\infty\leq \Vert f_\epsilon|_\La\Vert_\infty/\cos(\rho+\delta)$. By letting $\epsilon\to 0$, we obtain (3). $\Box$
       
  \newpage    
    \medskip\noindent{\bf 3.
Proof of Lemma 1} 

\medskip

1. The proof in \cite{c} is based on the following result from \cite{d} (for some extension see \cite{h}):
{\sl Let $f\in B_{[-\tau,\tau]}$ be a real function satisfying $-1\leq f(x)\leq 1$ for all $x\in\R$.
Then for every real $a$ the function $\cos (\tau z + a) - f(z)$ vanishes identically or else it
has only real zeros.  Moreover it has a zero in every interval where $\cos (\tau z + a)$ varies between -1 and 1 and all the zeros are
simple, except perhaps at points on the real axis where $f(x) = \pm  1.$}

Sketch of proof. We may assume $a=0$ and $\tau=1$. Consider the function 
$$
 f_\epsilon(z):=(1-\epsilon)\frac{\sin (\epsilon z)}{\epsilon z}f((1-\epsilon)z).
$$
One may  check that  $f_\epsilon\in B_{[-1,1]}$, $-1<f(t)<1, t\in\R,$ and that the estimate holds
$$
|f_\epsilon (z)|\leq \frac{e^{ |y|}}{\epsilon |z|}, z=x+iy \in\C.
$$ This shows that $|f_\epsilon(z)|<|\cos z|$ when $z$ lies on a rectangular contour $\gamma$ consisting of segments of the lines $x = \pm N\pi, y= \pm N,$ where $N$ is every large enough integer. By Rouch\`e's theorem, the function $\cos z-f_\epsilon(z)$
has the same number of zeros in $\gamma$ as $\cos  z$, that is, $2N$ zeros. On the
real axis $|f_\epsilon|\leq 1-\epsilon $. Hence,  $\cos z-f_\epsilon(z)$ is alternately plus and
minus at the $2N+1$ points $k\pi$, $|k|\leq N,$ so  it has  $2N$ real zeros inside $\gamma$. Taking larger values of $N$ we see that $\cos z-f_\epsilon(z)$ has exclusively real and simple
zeros, which lie in the intervals $(k\pi,(k+1)\pi)$.

The zeros of $\cos z-f(z)$ are limit points of the zeros of $\cos z-f_\epsilon(z)$  as $\epsilon\to 0$. Thus $\cos  z-f(z)$ cannot
have non-real zeros. Moreover, it has an infinite number of real zeros
which are all simple, except those at the points $k\pi$ iff $f(k\pi) = ( - 1)^{k}.$
Every interval $k\pi<z<(k + 1)\pi$ at the endpoints of which
$| f (t) | < 1$ contains exactly one zero. If $f(k\pi) = (-1)^{k}$, we have a
double zero at $k\pi$ but no further zeros in the interior or at the endpoints
of the interval $((k - 1)\pi, (k+ 1)\pi). $

2. It suffices to prove Lemma 1 for real functions $f\in B_{[-\tau,\tau]}$. Since $f$ has a local maximum at $t=0$, the function $f(t)-\cos\tau t$ has a repeated zero at $t=0. $
By the discussion above we see that either $f(t)$ is identically equal to $\cos\tau t$ or $f(t)-\cos\tau t$ does not vanish on $[-\pi/\tau,0)\cup(0,\pi/\tau]$. Since $|f(\pi)|\leq 1$, it follows that  $f(t)>\cos\tau t$ on each of the intervals $[-\pi/\tau,0)$ and $(0, \pi/\tau]$.
$\Box$

\end{large} 

\end{document}